\documentclass{amsart}
\usepackage{amsmath, amssymb, amsthm, graphicx}

\theoremstyle{plain}
\newtheorem{prop}{Proposition}
\newtheorem{coro}[prop]{Corollary}
\newtheorem{lemm}[prop]{Lemma}

\def\Reff#1; #2; #3; #4; #5; #6; #7\par{%
\bibitem{#1} #2, {\it #3}, #4 {\bf #5} (#6) #7}

\def\Ref#1; #2; #3; #4\par{%
\bibitem{#1} #2, {\it #3}, #4}

\def\ge{\geqslant}

\def\inv{^{-1}}
\def\lcm{{\rm LCM}}

\let\pp=\ldots

\begin{document}

\author{Patrick DEHORNOY}
\address{Laboratoire de Math\'ematiques Nicolas
Oresme UMR 6139\\ Universit\'e de Caen,
14032~Caen, France}
\email{dehornoy@math.unicaen.fr}
\urladdr{//www.math.unicaen.fr/\!\hbox{$\sim$}dehornoy}

\title{The group of fractions of a torsion free lcm monoid is torsion free}



\begin{abstract}
We improve and shorten the argument given in~\cite{Dfz}
(this journal, vol.~210 (1998) pp~291--297). In particular,
the fact that Artin braid groups are torsion free
now follows from Garside's results almost immediately.
\end{abstract}

\maketitle

An algebraic proof of the fact that Artin braid groups are
torsion free was given in~\cite{Dfz}. The aim of this note is
to observe that the proof of~\cite{Dfz}, which uses group
presentations and words, is unnecessarily complicated and
requires needless hypotheses. A shorter and better
argument can be given that uses elementary properties of
the lcm operation only.

For~$x, y$ in a monoid~$M$, we say that $y$ is a {\it right
multiple} of~$x$ if $y = x z$ holds for some~$z$
in~$M$; we say that $z$ is a right least common multiple,
or {\it right lcm}, of~$x$ and~$y$ if $z$ is a right multiple
of~$x$ and~$y$ and any common right multiple
of~$x$ and~$y$ is a right multiple of~$z$. The result
we prove here is:

\begin{prop}\label{P:Main}
Assume that $G$ is a group and $M$ is a
submonoid of~$G$ such that $M$ generates~$G$
and, in~$M$, any two elements admit a right lcm.
Then the torsion elements of~$G$ are the
elements~$x t x\inv$ with $x$ in~$M$
and $t$ a torsion element of~$M$.
\end{prop}

\begin{coro}\label{C:Tors}
Under the previous hypotheses, $G$ is torsion free if
and only if $M$ is torsion free. In particular,
a sufficient condition for $G$ to be torsion free
is that $M$ contains no invertible element but~$1$.
\end{coro}

In comparison with~\cite{Dfz}, the current result eliminates
a useless Noetherianity hypothesis. In this way, the result
directly extends the standard result that a left-orderable
group is torsion free. Indeed, if $<$ is a linear ordering
on~$G$ that is compatible with left multiplication, the
submonoid defined by $M = \{x \in G; x \ge 1\}$ is eligible
for Corollary~\ref{C:Tors}, as, in~$M$, the
element~$\sup(x, y)$ is a right lcm of~$x$ and~$y$. 

The above results apply to Artin's braid
group~$B_n$, as, according to Garside's theory, the
submonoid~$B_n^+$ admits unique right lcm's.
Alternatively, one could also use the dual monoid
of~\cite{BKL}, or some more exotic monoids.
Artin--Tits groups of finite Coxeter type, and, more
generally, all Garside groups \cite{Dgk} are eligible,
as well as, for instance, Richard Thompson's
group~$F$~\cite{CFP}. All lattice-ordered groups
of~\cite{GlC} also are eligible, but, then, the result is
trivial: the point is that, here, we only assume
one-sided compatibility between multiplication and
ordering.

In order to prove Proposition~\ref{P:Main}, we begin with
two simple observations about lcm's. First, a right lcm need
not be unique, but the set, here denoted~$\lcm(x, y)$, of all
right lcm's of two elements~$x, y$ is easily described:

\begin{lemm}\label{L:Uniq}
Assume that $M$ is a left cancellative monoid, and $z$ is a
right lcm of two elements~$x, y$ of~$M$. Then $\lcm(x,
y)$ consists of all elements of the form~$zu$ with~$u$
an invertible element of~$M$.
\end{lemm}

\begin{proof}
If~$u$ is invertible in~$M$, the element~$zu$ is a
right multiple of~$x$ and~$y$, and $z$ is a right
multiple of~$zu$, so $zu$ is a right lcm of~$x$
and~$y$. Conversely, let
$z'$ be an arbitrary right lcm of~$x$ and~$y$.
There must exist~$u, u'$ satisfying $z' = zu$ and
$z = z'u'$, hence $z = z u u'$ and $z' = z'u' u$, and
we deduce $u u' = u' u = 1$.
\end{proof}

\begin{lemm}\label{L:Rlcm}
Assume that $M$ is a left cancellative monoid, and that,
in~$M$, we have $x y'_1 = y_1 x' \in \lcm(x, y_1)$
and $x' y'_2 = y_2 x'' \in \lcm(x', y_2)$. Then we
have $x y'_1 y'_2 = y_1 y_2 x'' \in \lcm(x, y_1
y_2)$.
\end{lemm}

\begin{proof}
First  we have $x y'_1 y'_2 = y_1 x' y'_2 = y_1 y_2 x''$, so
this element is a common right multiple of~$x$ and~$y_1
y_2$. Assume that $z$ is a right multiple of~$x$
and of~$y_1 y_2$, say $z = y_1 y_2 z'$. Then $z$ is
a right multiple of~$x$ and~$y_1$, hence of~$y_1 x'$,
say $z = y_1 x' z''$. Cancelling~$y_1$ on the left, we
obtain $y_2 z' = x' z''$, so $y_2 z'$ is a right multiple
of~$y_2 x''$, and $z$ is a right multiple of~$y_1 y_2 x''$.
\end{proof}

\begin{proof}[Proof of Proposition~\ref{P:Main}]
The condition is obviously sufficient and the only
problem is to prove that it is necessary. As $M$ is a
submonoid of a group, it admits cancellation, and, as
any two elements of~$M$ admit a common right
multiple, $M$ satisfies the Ore conditions on the
right, and $G$ is a group of right fractions of~$M$.
Let $z$ be an arbitrary element of~$G$. Write $z =
x_1 y_1\inv$ with $x_1, y_1$ in~$M$, and,
inductively, choose
$x_2, y_2, x_3, y_3, \pp$ in~$M$ satisfying $x_i
y_{i+1} = y_i x_{i+1} \in \lcm(x_i, y_i)$.
We claim that, for all positive~$k, \ell$, we have
\begin{gather}
\label{E:1}
x_1 \pp x_k y_{k + 1} \pp y_{k+\ell} = 
y_1 \pp y_\ell x_{\ell + 1} \pp x_{\ell+k}
\in \lcm(x_1 \pp x_k, y_1 \pp y_\ell),\\
\label{E:2}
z = (x_1 \pp x_k) (x_{k+1} y_{k+1}\inv)
(x_1 \pp x_k)\inv,\\
\label{E:3}
z^k = (x_1 \pp x_k) (y_1 \pp y_k)\inv.
\end{gather}
Indeed, \eqref{E:1} follows from
Lemma~\ref{L:Rlcm} inductively; for~\eqref{E:2}
and~\eqref{E:3}, for each~$i$, we have
$y_i\inv x_i = x_{i+1} y_{i+1}\inv$
by construction, and we deduce
\begin{align*}
zx_1 \pp x_k 
&= (x_1 y_1\inv) x_1 \pp x_k
= x_1 (x_2 y_2\inv) x_2 \pp x_k = \pp
= x_1 \pp x_k( x_{k+1} y_{k+1}\inv),\\
z^k  y_1 \pp y_k 
&= (x_1 y_1\inv )^k  y_1 \pp y_k
 = x_1 (y_1\inv x_1)^{k-1}  y_2 \pp y_k
 = x_1 (x_2 y_2\inv)^{k-1}  y_2 \pp y_k \\
& = x_1 x_2 (y_2\inv x_2)^{k-2} y_3 \pp y_k 
= x_1 x_2 (x_3 y_3\inv)^{k-2} y_3 \pp y_k
= \pp = x_1 \pp x_k.
\end{align*}
Now assume $z^p = 1$, and let $t = x_{p+1}
y_{p+1}\inv$. By~\eqref{E:2}, we have $z = x t
x\inv$ with $x = x_1 \pp x_p \in M$.
Relation~\eqref{E:3} implies
\begin{equation}\label{E:4}
x_1 \pp x_p = 
y_1 \pp y_p
\in \lcm(x_1 \pp x_p, y_1 \pp y_p).
\end{equation}
Comparing Relations~\eqref{E:1}---with $k = \ell =
p$---and~\eqref{E:4}, we deduce from
Lemma~\ref{L:Uniq} that $y_{p+1}
\pp y_{2p}$, hence $y_{p+1}$ as well, is invertible
in~$M$. Therefore $t$ belongs to~$M$, and, as
$z$ and~$t$ are conjugates, $z^p = 1$ implies
$t^p = 1$.
\end{proof}

In lattice-ordered groups, the next result after
torsion freeness is that $x^p = y^p$ implies that $x$
and~$y$ are conjugate. This result does {\it not}
extend to our current framework: the group $\langle
x,y\,;\, x^2=y^2\rangle$ satisfies all hypotheses of
Proposition~1 but $x$ and~$y$ are not conjugate
there.

\end{document}